\documentclass[11pt,a4paper]{amsart}
\usepackage{amscd,amssymb,graphicx}
\usepackage{amsfonts,amsmath,amssymb,epsfig}
\usepackage{xypic}
\xyoption{all}
\input{xypic}
\newdir{ >}{{}*!/-8pt/\dir{>}}
\def\R{\mathbb{R}}
\def\N{\mathbb{N}}

\def\Z{\mathbb{Z}}

\def\fa{{\mathfrak a}}
\def\fb{{\mathfrak b}}
\def\fc{{\mathfrak c}}

\def\fg{{\mathfrak g}}
\def\fh{{\mathfrak h}}

\def\fk{{\mathfrak k}}

\def\fm{{\mathfrak m}}

\def\Strut{\rule[-1ex]{0ex}{5ex}}
\def\pr{\noindent $\bf{Proof.}$\quad}     
\def\fin{\hfill$\square$\\}           
\newtheorem{theo}{Theorem}

\newtheorem{prop}{Proposition}

\theoremstyle{remark}
\newtheorem{rem}{Remark}
\author{Alice Fialowski}
\address{E\"otv\"os Lor\'and University\\ Budapest, Hungary}
\email{fialowsk@cs.elte.hu}

\author{Friedrich Wagemann}
\address{Universit\'e de Nantes\\ Nantes, France}
\email{wagemann@math.univ-nantes.fr}

\begin{document}
\title[Associative algebra deformations of the Hopf algebra \mbox{$\mathcal{H}_1$}] 
{Associative algebra deformations of Connes-Moscovici's Hopf algebra
\mbox{$\mathcal{H}_1$}}

   


\begin{abstract}
We compute the second Hochschild cohomology space $HH^2(\mathcal{H}_1)$
of Connes-Moscovici's Hopf algebra $\mathcal{H}_1$, giving the infinitesimal 
deformations (up to equivalence) of the associative structure. The
space
$HH^2(\mathcal{H}_1)$ is shown to be one dimensional. 
\end{abstract}

\maketitle



\section*{Introduction}

Recently, there has been much interest in Connes-Moscovici's Hopf algebra
\mbox{$\mathcal{H}_1$} in relation to deformation quantization. Recall that
\mbox{$\mathcal{H}_1$} was constructed by Connes and Moscovici in \cite{ConMos}
in order to formalize the transverse symmetries of a codimension $1$ foliation.
As an associative algebra, \mbox{$\mathcal{H}_1$} is the universal enveloping 
algebra $U\fh$ of a certain Lie algebra $\fh$ closely related to the 
`$ax+b$'-group. In \cite{ConMos}, \cite{ConMos2}, the authors study the Hopf 
cyclic cohomology of (within others) \mbox{$\mathcal{H}_1$}, and associate
thereby characteristic classes to a codimension $1$ foliation. This
turns out to
be closely related to Gelfand-Fuchs cohomology classes. 

Notwithstanding their initial work \cite{ConMos}, Connes and Moscovici went on
remarking that Rankin-Cohen brackets, which were known to give associative 
deformations on spaces of modular forms, can give such deformations
on all algebras \mbox{$\mathcal{H}_1$} acts on \cite{ConMos3}. 
The answer is given in the framework 
of universal deformation formulae (UDF), see \cite{GiaZha}. 
The universality of this kind of deformation is
the fact that all \mbox{$\mathcal{H}_1$}-module algebras 
inherit a deformation of the multiplication.

In further work, Bieliavsky, Tang and Yao \cite{BieTanYao}
refind the Rankin-Cohen 
starproduct in the context of Fedosov deformation quantization and showed that 
the Rankin-Cohen UDF on \mbox{$\mathcal{H}_1$} 
is related to the Moyal-Weyl starproduct on the universal enveloping algebra 
of the `$ax+b$'-group. 
Still more recently, Tang and Yao \cite{TanYao}
showed that the \mbox{$\mathcal{H}_1$}-action 
does not have to be projective in order to define the Rankin-Cohen 
starproduct.

Our contribution to the subject is a computation of the space of infinitesimal
deformations of the associative product of \mbox{$\mathcal{H}_1$}, namely the 
Hochschild cohomology space
$HH^2(\mathcal{H}_1)$. It turns out to be $1$-dimensional, and this
shows the
unicity (up to equivalence) of the infinitesimal term of any associative 
formal deformation. 

The computation is performed using Lie algebra cohomology methods and spectral 
sequences. Indeed, Hochschild cohomology of a universal enveloping algebra 
boils down to Lie algebra cohomology. Here it is the Lie algebra cohomology 
space $H^2(\fh;\mathcal{H}_1^{\rm ad})$ of $\fh$ with values in 
$\mathcal{H}_1$ using the adjoint action. To go on with the computation, we 
use the link between $\fh$ and a Lie algebra called $\fm_0$. The
algebra $\fm_0$ is one 
of the three positively graded, infinite dimensional filiform Lie algebras 
(up to isomorphism) \cite{Fia1} and has been studied intensively in 
\cite{FialMill}, \cite{FiaWag}. In fact, $\fm_0$ is an ideal 
of $\fh$, and thus we may use the Hochschild-Serre spectral sequence to compute
$H^2(\fh;\mathcal{H}_1^{\rm ad})$ from the various 
$H^p(\fh\,/\,\fm_0;H^q(\fm_0;\mathcal{H}_1^{\rm ad}))$ with $p+q=2$. 

These latter spaces $H^p(\fh\,/\,\fm_0;H^q(\fm_0;\mathcal{H}_1^{\rm ad}))$ 
are rather easily deduced from the knowledge of the spaces
$H^q(\fm_0;\mathcal{H}_1^{\rm ad})$, $q=0,1,2$, which are computed using the 
Feigin-Fuchs spectral sequence - for an introduction to it, see \cite{MilWag}.

Let us give a detailed account on the content of this paper: in Section
\ref{Connes},
we define $\mathcal{H}_1$ (which we call in the following only $\mathcal{H}$)
in Connes and Moscovici's original context. In $1.2$, $\fh$ is defined by 
generators and relations, and $1.3$ gives some background material on $\fm_0$.
Section \ref{link} gives the link between $\fh$ and $\fm_0$, commenting on other
points of view which permit, for example, to associate a pro-Lie group to $\fh$.
Section \ref{outline} gives an outline of the cohomology computations: in 
Section \ref{Hoch}
the link between Hochschild cohomology of $U\fg$ and Chevalley-Eilenberg 
cohomology of $\fg$ is recalled. Section \ref{spectral} recalls the Hochschild-Serre
spectral sequence, and Section \ref{Feigin} the Feigin-Fuchs spectral sequence.
Section \ref{Comp} treats the computations: in Section \ref{cohom}, we compute the spaces
$H^q(\fm_0;\mathcal{H}_1^{\rm ad})$, $q=0,1,2$ corresponding to Propositions 
\ref{zerocoh}, \ref{onecoh}, and \ref{twocoh}. In Section \ref{ad}, we deduce then 
$H^p(\fh\,/\,\fm_0;H^q(\fm_0;\mathcal{H}_1^{\rm ad}))$ with $p+q=2$ mainly
by degree arguments. Observe here that while $\fh$ has a basis consisting
of eigenvectors with respect to $Y\in\fh$, the grading of $\fm_0$ is not inner:
that is the fundamental difference between $\fh$ and $\fm_0$ which renders 
computations easy for $\fh$ and difficult for $\fm_0$. It is the Feigin-Fuchs
spectral sequence which helps out.

The main theorems of this paper are Theorems \ref{coh1} and \ref{coh2}
in Section \ref{ad} which state the result of the cohomology
computations of $HH^1(\mathcal{H})$ and $HH^2(\mathcal{H})$ in terms of
Lie algebra cocycles whose classes  generate these spaces.

Our approach is new in the sense that up to now, only the deformation theory of
algebras $A$ where $\mathcal{H}$ acts on (such that 
$h(ab)=m(\triangle(h)(a\otimes b))$ holds for all $h\in\mathcal{H}$ and all 
$a,b\in A$) has been regarded. This is the first step towards a deformation 
theory of the algebra $\mathcal{H}$ itself. In this context,
it would be interesting to construct a starproduct whose infinitesimal term 
represents the generator of $HH^2(\mathcal{H})$. Such a starproduct would then
be automatically the miniversal deformation of $\mathcal{H}$.\\

\noindent{\bf Acknowledgements:}\quad FW thanks Dmitri Millionschikov for 
explaining the Feigin-Fuchs spectral sequence to him. Both authors
thank the R\'enyi
Institute in Budapest where a first version of these results was
exposed. We also thank Xiang Tang for answering kindly our questions.

\section{Preliminairies on \mbox{$\mathcal{H}_1$} and ${\mathfrak m}_0$}

\subsection{Connes-Moscovici's Lie and Hopf algebras}\label{Connes}

Let us recall some basics about Connes-Moscovici's Hopf algebra
\mbox{$\mathcal{H}_1$}. It was introduced in \cite{ConMos} in the study of 
the transverse structure of foliations. 

Let $M$ be an $n$-dimensional smooth oriented manifold, $\nabla$ be a flat 
affine connection on $M$, and $\Gamma$ be a pseudogroup of local 
diffeomorphisms on $M$, respecting the orientation. Such a pseudogroup arises
for example in the presence of an oriented foliation. 

Denote by $F^+$ the oriented frame bundle on $M$. Define \cite{ConMos}
$${\mathcal A}:={\mathcal C}^{\infty}_c(F^+)\rtimes\Gamma,$$
the crossed product of the algebra ${\mathcal C}^{\infty}_c(F^+)$ of smooth 
sections of $F^+$ with compact support with the pseudogroup $\Gamma$. 
 
There are three kinds of generators of ${\mathcal A}$:

\noindent(1)\quad Generators $Y^i_l$:

$Gl^+(n,\R)$ acts on $F^+$, and this action, as it commutes with that of 
$\Gamma$, extends to ${\mathcal A}$. The infinitesimal generators (or 
fundamental vector fields) of this action are by definition the $Y^i_l$,
$i,l=1,\ldots,n$.    

\noindent(2)\quad Generators $X_i$:

The flat connection $\nabla$ permits to lift vector fields which are tangent 
to $M$, to $F^+$: we therefore get horizontal vector fields $X_i$ on $F^+$,
$i=1,\ldots,n$.

\noindent(3)\quad Generators $\delta^t_{rs,i_1,\ldots,i_l}$:

In case the local diffeomorphisms are affine, the $X_i$-s commute with the action
of elements of $\Gamma$, but in general, they do not. One finds \cite{ConMos}
$$X_i(ab)\,=\,X_i(a)b+aX_i(b)+\delta^k_{ij}Y^j_k(b),$$
for all $a,b\in{\mathcal A}$ and the vertical fields $Y^j_k$, infinitesimal 
generators of the action of $Gl^+(n,\R)$ on $F^+$. Iterated brackets of
the $X_i$ and the $\delta^t_{rs}$ yield 
$$\delta^t_{rs,i_1,\ldots,i_l}\,:=\,[X_{i_1},\ldots,[X_{i_l},\delta^t_{rs}]
\ldots].$$

The space generated by the $Y^i_l$, $i,l=1,\ldots,n$, $X_i$, $i=1,\ldots,n$,
and $\delta^t_{rs,i_1,\ldots,i_l}$, $r,s,t=1,\ldots,n$, $l\in\N$, is closed 
under the Lie bracket (by construction), and yields therefore a Lie algebra
${\mathfrak h}(n)$. Its enveloping algebra is denoted by \mbox{$\mathcal{H}_n$}.

Connes and Moscovici \cite{ConMos} endow the associative algebra
\mbox{$\mathcal{H}_n$} with a coproduct in such a way that it acts on the 
algebra ${\mathcal A}$:
$$h(ab)\,=\,\sum h_0(a)h_1(b),$$
(in Sweedler notation for the coproduct) for all $a,b\in{\mathcal A}$ and all
$h\in\mathcal{H}_n$. They then construct an antipode and \mbox{$\mathcal{H}_n$}
becomes a Hopf algebra; we call it the {\it Connes-Moscovici Hopf algebra}
\mbox{$\mathcal{H}_n$}.

We will be only concerned with the $1$-dimensional case, i.e. $n=1$, in which
case we will write more simply ${\mathfrak h}$ and \mbox{$\mathcal{H}$} for
Connes-Moscovici's Lie and Hopf algebras.

We will not recall the significance of \mbox{$\mathcal{H}_n$} in foliation 
theory, as this would lead us too far afield. 

\subsection{Connes-Moscovici's Lie algebra ${\mathfrak h}$}

Algebraically speaking, ${\mathfrak h}$ is a Lie algebra\footnote{Because of 
their geometric origin, we take Lie algebras over $\R$ here; the algebraic 
results are of course valid over much more general fields.} generated by
the countably infinite set of generators $X$, $Y$ and $\delta_i$, 
$i=1,2,\ldots$, with the relations
$$\begin{array}{rclrclrcl}
[Y,X]&=&X, & [Y,\delta_r]&=&r\delta_r, & [X,\delta_r]&=&\delta_{r+1},
\end{array}$$
for all $r=1,2,\ldots$. By convention, we only write the non-trivial 
relations, i.e. all brackets which are not displayed, are zero. In particular,
one has $[\delta_r,\delta_s]=0$ for all $r,s=1,2,\ldots$. 

Observe that $\fh$ is a graded Lie algebra: $Y$ has degree
$0$, $X$ and $\delta_1$ have degree $1$, and $\delta_i$ has degree $i$
for $i\geq 1$.

\subsection{The infinite dimensional filiform Lie algebra ${\mathfrak m}_0$}

The significance of ${\mathfrak m}_0$ resides in the classification of 
infinite dimensional $\N$-graded Lie algebras 
$\fg=\bigoplus_{i=1}^{\infty}\fg_i$ with one-dimensional homogeneous components
$\fg_i$ and two generators (over a field of characteristic zero)
such that $[\fg_1,\fg_i]=\fg_{i+1}$. A. Fialowski 
showed in \cite{Fia1} that any Lie algebra of this type must be isomorphic to 
$\fm_0$, $\fm_2$ or $L_1$. We call these Lie algebras infinite dimensional 
filiform Lie algebras in analogy with the finite dimensional case where the 
name was coined by M. Vergne in \cite{Ver}. Here $\fm_0$ is given by generators
$e_i$, $i\geq 1$, and relations $[e_1,e_i]=e_{i+1}$ for all $i\geq 2$, $\fm_2$ 
with the same generators by relations $[e_1,e_i]=e_{i+1}$ for all $i\geq 2$, 
$[e_2,e_j]=e_{j+2}$ for all $j\geq 3$, and $L_1$ with the same generators is 
given by the relations $[e_i,e_j]=(j-i)e_{i+j}$ for all $i,j\geq 1$. $L_1$
appears as the positive part of the Witt algebra given by generators $e_i$ for
$i\in\Z$ with the same relations $[e_i,e_j]=(j-i)e_{i+j}$ for all $i,j\in\Z$.

What matters most for the present discussion is the cohomology with trivial 
coefficients of $\fm_0$. It has been computed in \cite{FialMill} and will be
recalled in Theorem \ref{bigraded}.

\subsection{The link between $\fh$ and $\fm_0$} \label{link}

Two ways of viewing $\fh$  are of interest: first, define a map 
$i:\fm_0\to\fh$ by sending $e_1$ to $X$ and $e_i$ to $\delta_i$ for $i\geq 2$.
A short inspection of the relation shows that via $i$, $\fm_0$ becomes a 
subalgebra of $\fh$, and one has a short exact sequence:
$$0\to\fm_0\stackrel{i}{\to}\fh\to{\rm span}(Y,\delta_1)\to 0.$$
This determines $\fh$ as a general extension of a 
$2$-dimensional Lie algebra $\fa:={\rm span}(Y,\delta_1)$ (with the only 
non-trivial relation $[Y,\delta_1]=\delta_1$) by $\fm_0$. The
term "general extension" means here that it is neither a central, nor an 
abelian extension. The algebra $\fa$ is the Lie algebra of the `$ax+b$'-group, see
\cite{BieTanYao}.

The second way of viewing $\fh$ is as a trivial abelian extension of the 
$2$-dimensional Lie algebra $\fb:={\rm span}(X,Y)$ (with the relation 
$[Y,X]=X$) by the infinite dimensional abelian Lie algebra 
$\fc:={\rm span}(\delta_i\,|\,i\geq 1)$:
$$0\to\fc\to\fh\to\fb\to 0.$$

In the following, we will exploit the first point of view to compute some
cohomology of $\fh$; let us remark here that the second point of view makes
it possible to associate an infinite dimensional Lie group to the Lie algebra 
$\fh$. Indeed, truncating the infinite dimensional $\fc$ to a finite 
dimensional $\fc_k:={\rm span}(\delta_i\,|\,i=1,\ldots,k)$ defines a Lie
algebra $\fh_k$ by
$$0\to\fc_k\to\fh_k\to\fb\to 0.$$
Then $\fh$ is obviously a projective limit of the finite dimensional
Lie algebras $\fc_k$, and thus there is a Lie group associated to $\fh$
by the following  theorem \cite{HofMor} (Lie's third theorem for
pro-Lie groups):

\begin{theo}[Hofmann-Morris]
Given a profinite Lie algebra, i.e. a projective limit of finite dimensional
Lie algebras, there is a connected profinite (possibly infinite dimensional)
Lie group whose Lie algebra is the given one.
\end{theo}

It would be interesting to exploit this last remark in order to establish
some relationship between the pro-Lie group associated to $\fh$ and 
\mbox{$\mathcal{H}$}.
\bigbreak

\section{Outline of the cohomology computation}\label{outline}

The aim is to compute the second Hochschild cohomology 
$HH^2(\mathcal{H})$ of the associative algebra \mbox{$\mathcal{H}$}
in order to determine the different infinitesimal deformations of 
\mbox{$\mathcal{H}$} as an associative algebra.
We will do this in three steps. 

\subsection{From Hochschild to Chevalley-Eilenberg}\label{Hoch}

Recall that {\it Hochschild cohomology} of an associative algebra $A$ over a 
field $k$ with values in an $A$-bimodule $M$ is just
$$HH^*(A;M)\,:=\,Ext_{A\otimes A^{\rm opp}}^*(k,M).$$
In case $M=A$ with the usual bimodule structure, we will write $HH^*(A)$ 
instead of $HH^*(A;A)$. The second cohomology group $HH^2(A)$ 
classifies infinitesimal deformations of the algebra structure of $A$
\cite{Ger}.

Recall also that the {\it Chevalley-Eilenberg cohomology} of a Lie
algebra $\fg$ over $k$ with values in a $\fg$-module $N$ is by definition
$$H^*(\fg;N)\,:=\,Ext_{U\fg}^*(k,N),$$
where $U\fg$ is the universal enveloping algebra of $\fg$. 

Let $M$ be a $U\fg$-bimodule, and denote by $M^{\rm ad}$ the right 
$\fg$-module defined by
$$m\cdot x\,=\,mx - xm.$$
It is explained in the book of Cartan-Eilenberg \cite{CarEil} that the 
change-of-rings functor associated to the map 
$\fg\to U\fg\otimes U\fg^{\rm opp}$ given by
$x\mapsto x\otimes 1 - 1\otimes x$ leads to an isomorphism
$$HH^*(U\fg;M)\,\cong\,H^*(\fg;M^{\rm ad})$$
in what they call the inverse process. While the homomorphism inducing this 
isomorphism is well-understood in one direction, its inverse does not appear 
in the literature. Indeed, the isomorphism is induced by the antisymmetrization
map
$$x_1\wedge\ldots\wedge x_p\mapsto\sum_{\sigma\in S_p}(-1)^{{\rm sgn}(\sigma)}
x_{\sigma(1)}\otimes\ldots\otimes x_{\sigma(p)}$$
between standard resolutions, where $x_1,\ldots,x_p\in\fg$. 

In this way, in order to compute $HH^2(\mathcal{H})$ for the Connes-Moscovici 
Hopf algebra $\mathcal{H}$, it is enough to compute the Lie algebra cohomology
space $H^2(\fh;\mathcal{H}^{\rm ad})$. 
 
\subsection{From $\fh$ to $\fm_0$ via the Hochschild-Serre spectral sequence}
\label{spectral}

In the second step, we will use the short exact sequence
$$0\to\fm_0\to\fh\to{\rm span}(Y,\delta_0)\to 0$$
in order to reduce the computation to one for $\fm_0$. This is done by the 
Hochschild-Serre spectral sequence: given an ideal $\fk\subset\fg$ in a Lie
algebra $\fg$ and a $\fg$-module $L$, there is a filtration on the space of 
cochains $C^*(\fg;L)$ which induces a spectral sequence with
$$E_2^{p,q}\,=\,H^p(\fg/\fk;H^q(\fk;L))$$
converging to 
$$E_{\infty}^{p,q}\,=\,H^{p+q}(\fg;L),$$
see for example \cite{Fuchs}. 

In our case, we therefore have to compute the spaces
$H^q(\fm_0;\mathcal{H}^{\rm ad})$ for $q=0,1,2$, and then the cohomology of 
$\fa:={\rm span}(Y,\delta_0)$ with values in these spaces in order to 
determine the $E_2$-term. As the latter
computation is rather easy, we are left with computing 
$H^q(\fm_0;\mathcal{H}^{\rm ad})$ which will be done in the third step. 

\subsection{The computation for $\fm_0$ via the Feigin-Fuchs spectral sequence}
\label{Feigin}

$H^q(\fm_0;\mathcal{H}^{\rm ad})$ will be computed using the Feigin-Fuchs
spectral sequence. This is a tool which is available only for $\N$-graded
Lie algebras with values in a non-negatively graded module. 
It has been introduced 
in \cite{FeiFuc}.
The preprint \cite{MilWag} is meant to be a readable introduction to this 
subject.

Let $\fk$ be an $\N$-graded Lie algebra, i.e. 
$\fk=\bigoplus_{i=1}^{\infty}\fk_i$ with $[\fk_i,\fk_j]\subset \fk_{i+j}$,
and let $L$ be a non-negatively graded $\fk$-module, i.e. 
$L=\bigoplus_{i=0}^{\infty}L_i$ with $\fg_i\cdot L_j\subset L_{i+j}$.
Using the gradation, one introduces a filtration on $C^*(\fk;L)$ which induces 
a spectral sequence with
$$E_1^{p,q}\,=\,H^{p+q}(\fk)\otimes L_p$$
converging to the cohomology $H^*(\fk;L)$.   

The strategy of computation in this third step is thus clear: using the known 
results on the cohomology with trivial coefficients of $\fm_0$ from
\cite{FialMill}, one has to follow them through the Feigin-Fuchs spectral 
sequence. 
\bigbreak

\section{Computation and Results}\label{Comp}

Obviously, we will begin with the third step of the outline.

\subsection{Computation of $H^q(\fm_0;\mathcal{H}^{\rm ad})$}
\label{cohom}

The Feigin-Fuchs spectral sequence will become important for
$H^q(\fm_0;\mathcal{H}^{\rm ad})$ for $q=1,2$, while we first compute 
$H^0(\fh;\mathcal{H}^{\rm ad})$ by elementary methods.

As it is well known, $H^0(\fh;\mathcal{H}^{\rm ad})$ is the subspace of
$\fh$-invariants of $\mathcal{H}^{\rm ad}$. We have:

\begin{prop}\label{zerocoh}\qquad $\displaystyle H^0(\fh;\mathcal{H}^{\rm
ad})\,=\,\R$.
\end{prop}

\pr As $\fh$ is a graded Lie algebra and $\mathcal{H}^{\rm ad}$ is a graded 
$\fh$-module, and we may reason degree by degree.

In degree $0$, we have polynomials in $Y$ and $\R$, and it is then clear by
$[Y,X]=X$ that the invariants in degree $0$ are $\R$. 
In degree $1$, we have $X$ and $\delta_1$.
As $[Y,X]=X$ and $[X,\delta_1]=\delta_2$, invariants are $0$ in degree $1$. 

Now let $m\in\mathcal{H}^{\rm ad}$ be a homogeneous element of degree $n>1$
with $m\cdot a=0$ for all $a\in\fh$. Suppose a non-zero $\delta_i$ occurs in 
$m$. By the Poincar\'e-Birkhoff-Witt theorem, there is a maximal $k$ and 
a corresponding maximal $r=r(k)$ such that $m=m'\delta^r_k$. We have
\begin{eqnarray*}
\delta^r_k\cdot X&=&\delta^r_kX-X\delta^r_k\,=\,\delta^{r-1}_kX\delta_k - 
\delta^{r-1}_k\delta_{k+1}-X\delta^r_k\\
&=&\ldots = - r\delta_k^{r-1}\delta_{k+1},
\end{eqnarray*}
because the $\delta_i$-s commute with each other. Then $m\cdot X=0$ implies thus
$rm'\delta^{r-1}_k\delta_{k+1}=0$ with $r\geq 1$. One deduces $m'=0$, and 
$m=0$.

The last case is the one where there is no $\delta_i$ in $m$. Suppose
therefore $m=\sum_kY^kX^n$ ($m$ is supposed to be of degree $n>1$ !), and again
$m\cdot X=0$. This gives 
$$m\cdot X\,=\,\sum_km_kY^kX^{n+1}-\sum_km_kXY^kX^n\,=\,0,$$
and we want to commute $X$ to the right. We have the (binomial) formula:
$$XY^k\,=\,(Y^k-kY^{k-1}\pm \ldots +(-1)^k)X,$$
which is easily shown by induction. One deduces
$$\sum_km_k(-ky^{k-1}+\left(\begin{array}{c} k \\ 2 \end{array}\right)Y^{k-2}
\mp\ldots+(-1)^k)\,=\,0.$$
Looking at the highest power of $Y$, say $k_{\rm max}$, we get 
$m_{k_{\rm max}}=0$, and then $m_{k_{\rm max}-1}=0$ and so on. Finally, $m=0$.
\fin

Let us now compute $H^1(\fm_0;\mathcal{H}^{\rm ad})$ via the Feigin-Fuchs 
spectral sequence. For an introduction to this spectral sequence, see
\cite{MilWag}. 

First recall the cohomology $H^*(\fm_0)$ of $\fm_0$ with trivial coefficients,
see Theorem {\bf 3.4} in \cite{FialMill}. The generators of $\fm_0$ are still
denoted $e_i$, $i\geq 1$, and a dual ``basis'' is given by the $1$-cochains
$e^i$, $i\geq 1$. Fialowski and Millionschikov  show in \cite{FialMill}:

\begin{theo}\label{bigraded}
The bigraded cohomology algebra $H^*(\fm_0)=\bigoplus_{k,q}H^q_k(\fm_0)$ is 
spanned by the cohomology classes of the following homogeneous cocycles:
$$e^1,e^2,\omega(e^{i_1}\wedge\ldots\wedge e^{i_q}\wedge e^{i_q+1})=
\sum_{l\geq 0}(-1)^l({\rm ad}_{e_1}^*)^l(e^{i_1}\wedge\ldots\wedge e^{i_q})
\wedge e^{i_q+1+l},$$
where $q\geq 1$, $2\leq i_1<\ldots<i_q$.
\end{theo}

 Fialowski and Millionschikov describe also the dimension of the cohomology 
spaces and the multiplicative structure in detail.

Coming back to $H^1(\fm_0;\mathcal{H}^{\rm ad})$, we have (we do not 
distinguish cohomology classes and cocycles generating it):

\begin{prop}\label{onecoh}\qquad $\displaystyle
H^1(\fm_0;\mathcal{H}^{\rm ad})\,\cong\,H^1(\fm_0)\oplus 
\R(\delta_1\otimes e^1)$.\\
More precisely, in the Feigin-Fuchs spectral sequence, we have:
$$E_{\infty}^{1,0}\,\cong\,H^1(\fm_0)$$
and 
$$E_{\infty}^{0,1}\,=\,\R(\delta_1\otimes e^1).$$
\end{prop}

\pr The Feigin-Fuchs spectral sequence uses the graded structure of Lie algebra 
and module. The action by Lie algebra elements sends module elements 
necessarily in strictly upper degrees, so the first step, leading to 
$E_1^{p,q}$, is to exclude
action terms in the Lie algebra differential. It remains the differential
with trivial coefficients and the $E_1$-term in the Feigin-Fuchs spectral 
sequence is
$$E_1^{p,q}\,=\,(\mathcal{H}^{\rm ad})_q\otimes H^{p+q}(\fm_0).$$
The space $H^{1}(\fm_0)$ is generated by (the cohomology classes represented by) $e^1$ 
and $e^2$. Therefore
$$E_1^{0,1}\,=\,(\bigoplus_{n\geq 0}\R Y^nX\oplus\bigoplus_{n\geq 0}\R
Y^n\delta_1)\otimes_{\R}(\R e^1\oplus \R e^2),$$
and
$$E_1^{1,0}\,=\,(\R 1\oplus\bigoplus_{n\geq 1}\R Y^n)\otimes_{\R}
(\R e^1\oplus \R e^2).$$
The second step, leading to the differential $d_1^{p,q}$, is then to admit 
action terms with elements of degree $1$. Applied to a general cochain
$$c\,=\,\sum_{n\geq 0}a_n Y^nX\otimes e^1+\sum_{n\geq 0}b_n Y^n\delta_1\otimes 
e^1+\sum_{n\geq 0}c_n Y^nX\otimes e^2+\sum_{n\geq 0}d_n Y^n\delta_1\otimes 
e^2,$$
the only non-zero term is given by the action term involving $X$. One gets
with $i\geq 2$ (remember that $\fm_0$ is identified with 
$\R X\oplus\bigoplus_{i\geq 2}\R \delta_i$) 
$$d_1^{0,1}c(X,\delta_i)\,=\,-c(\delta_i)\cdot X,$$
because the other action term $c(X)\cdot \delta_i$ does not take values in 
$(\mathcal{H}^{\rm ad})_{1}$. By the special form of the cochain $c$, we must 
have $i=2$. Using $Y^nX\cdot X =  
\sum_{l=0}^{n-1}Y^lXY^{n-l-1}X$ and $Y^n\delta_1\cdot X = 
\sum_{l=0}^{n-1}Y^lXY^{n-l-1}\delta_1 - Y^n\delta_2$ (which are easily shown 
by induction), we get
$$d_1^{0,1}c(X,\delta_2)\,=\,\sum_{n\geq 1}c_n \sum_{l=0}^{n-1}Y^lXY^{n-l-1}X
-\sum_{n\geq 0}d_n \left(\sum_{l=0}^{n-1}Y^lXY^{n-l-1}\delta_1 - 
Y^n\delta_2\right).$$
Then the condition $d_1^{0,1}c(X,\delta_2)=0$ means that all $c_n$, $n\geq 1$, 
and all $d_r$, $r\geq 0$, must be zero. Indeed, the two sums must be finite, 
and starting with the highest power in $Y$, one shows one by one that all 
terms are zero. Therefore
$$E_2^{0,1}\,=\,\left(\bigoplus_{n\geq 1}\R Y^nX\oplus\bigoplus_{n\geq 0}\R
Y^n\delta_1\right)\otimes_{\R}\R e^1 \oplus\R X\otimes e^2,$$
because $dY=X\otimes e^1$ is a coboundary.

The next differential $d_2^{0,1}$ involves the action terms with elements of 
degree~$2$. We have to compute here the term 
$$d_2^{0,1}c(X,\delta_2)\,=\,c(X)\cdot \delta_2\,=\,\left(\sum_{n\geq 1}a_n 
Y^nX+\sum_{n\geq 0}b_n Y^n\delta_1\right)\cdot\delta_2+c_0Xe^2(X).$$
In the same way as before, using this time $Y^nX\cdot\delta_2=
2(1+Y+\ldots+Y^{n-1})\delta_2\delta_1$ and $Y^nX\cdot\delta_2=
2(1+Y+\ldots+Y^{n-1})\delta_2X + Y^n\delta_3$, we get that $b_0$ is arbitrary, 
while all $a_n$ for $n\geq 0$ and all $b_r$ for $r\geq 1$ are zero. Therefore
$$E_3^{0,1}\,=\,\R X\otimes e^2 \oplus\R \delta_1\otimes e^1.$$
A quick look at the area of non-zero terms in the first page of the spectral 
sequence shows that there are no non-zero terms strictly below the $q=0$ axis
and strictly to the left of the $p=-q$ antidiagonal:
$$
\begin{array}{lllll}
\Strut E_1^{-2,3} &  E_1^{-1,3} &  E_1^{0,3} &  E_1^{1,3} &  E_1^{2,3} \\
\Strut E_1^{-2,2}  &  E_1^{-1,2}  &  E_1^{0,2} &  E_1^{1,2} &  E_1^{2,2} \\
\Strut  E_1^{-2,1}=0 &  E_1^{-1,1} &  E_1^{0,1} &  E_1^{1,1} &  E_1^{2,1} \\
\Strut  E_1^{-2,0}=0 &  E_1^{-1,0}=0  &  E_1^{0,0} &  E_1^{1,0} &  E_1^{2,0} \\
\Strut  E_1^{-2,-1}=0 \quad &  E_1^{-1,-1}=0 \quad  &  E_1^{0,-1}=0\quad  &  
 E_1^{1,-1}=0\quad  &  E_1^{2,-1}=0  \\
\end{array}
$$
Therefore, all images of $d_2^{p,q}:E_2^{p,q}\to E_2^{p-1,q+2}$, 
$d_3^{p,q}:E_3^{p,q}\to E_3^{p-2,q+3}$, and so on, in $E_r^{0,1}$ are zero, 
but one has an infinite
number of outgoing non-zero differentials. It is thus more convenient at this 
stage to compute which combinations of $X\otimes e^2$ and 
$\delta_1\otimes e^1$ are actual cocycles. It turns out that $X\otimes e^2$ is 
not a cocycle (test on $e_2$ and $e_i$ with $i\geq 3$), but 
$\delta_1\otimes e^1$ is a cocycle. Finally 
$$E_{\infty}^{0,1}\,=\,\R \delta_1\otimes e^1.$$
 
Now let us compute $E_{\infty}^{1,0}$. We start from a cochain 
$$c\,=\,a 1\otimes e^1 + b 1\otimes e^2 + \sum_{n\geq 1}c_n Y^n\otimes e^1 + 
\sum_{n\geq 1}d_n Y^n\otimes e^2$$
in
$$E_1^{1,0}\,=\,(\R 1\oplus \bigoplus_{n\geq 1}\R Y^n)\otimes_{\R}
(\R e^1\oplus \R e^2).$$
As before, we have to compute $c(\delta_2)\cdot X$, and this gives here
$c(\delta_2)\cdot X=\sum_{n\geq 1}d_n Y^n\cdot X$. The first condition is thus 
$d_n=0$ for all $n\geq 1$. For $d_2^{1,0}$,
we have to compute $c(X)\cdot\delta_2$, which gives 
$c(X)\cdot\delta_2=\sum_{n\geq 1}c_n Y^n\cdot\delta_2$, and therefore the 
second condition is $c_n=0$ for all $n\geq 1$. 
It is clear that all combinations of $1\otimes e^1$ and 
$1\otimes e^2$ are actual cocycles, and we get therefore
$$E_{\infty}^{1,0}\,\cong\,H^1(\fm_0),$$
$H^1(\fm_0)$ being generated by $e^1$ and $e^2$.\fin

We go on computing $H^2(\fm_0;\mathcal{H}^{\rm ad})$ with the Feigin-Fuchs 
spectral sequence. 

\begin{prop}\label{twocoh}\qquad $\displaystyle
H^2(\fm_0;\mathcal{H}^{\rm ad})\,\cong\,H^2(\fm_0)$.\\
More precisely, in the Feigin-Fuchs spectral sequence, we have:
$$E_{\infty}^{2,0}\,\cong\,H^2(\fm_0)$$
and 
$$E_{\infty}^{1,1}\,=\,E_{\infty}^{0,2}\,=\,0.$$
\end{prop}

\pr The assertion $E_{\infty}^{2,0}\,\cong\,H^2(\fm_0)$ is clear in the 
same way that we computed $E_{\infty}^{1,0}\,\cong\,H^1(\fm_0)$: indeed,
the $E_1^{2,0}$-term is a sum of terms of the form 
$\R Y^n\otimes_{\R}H^2(\fm_0)$ for all $n$. The action with $X$ on some $Y^n$
gives once again a sum like $Y^n\cdot X = \sum_{l=0}^{n-1}Y^lXY^{n-l-1}$.
As the evaluation of the cochain on some elements from $\fm_0$ must be a 
finite sum, there is a term of highest degree in $Y$. By induction, we show 
as before that all terms involving $Y$ must be zero. The terms 
$\R 1\otimes_{\R}H^2(\fm_0)$ remain and give  
$$E_{\infty}^{2,0}\,\cong\,H^2(\fm_0).$$

Now let us look at $E_1^{1,1}\,=\,\mathcal{H}^{\rm ad}_1\otimes H^2(\fm_0)$
which may be written like
$$E_1^{1,1}\,=\,\left(\bigoplus_{n\geq 0}\R Y^nX\oplus\bigoplus_{n\geq 0}Y^n
\delta_1\right)\otimes_{\R}H^2(\fm_0),$$
and $H^2(\fm_0)$ is a countably infinite dimensional space with generators
$e^2\wedge e^3$, $e^2\wedge e^5 - e^3\wedge e^4$, and so on 
(cf \cite{FialMill}). Taking a cochain $c$ in $E_1^{1,1}$, it may be written 
like
\begin{eqnarray*}
c&=&\sum_{n\geq 0}a_n^1 Y^nX\otimes e^2\wedge e^3+\sum_{n\geq 0}b_n^1 
Y^n\delta_1\otimes e^2\wedge e^3  \\
&+&\sum_{n\geq 0}a_n^2 Y^nX\otimes (e^2\wedge e^5 - e^3\wedge e^4)+
\sum_{n\geq 0}b_n^2 Y^n\delta_1\otimes(e^2\wedge e^5 - e^3\wedge e^4)+\ldots.
\end{eqnarray*}
Then obviously 
$$d_1^{1,1}c(e_1,e_2,e_3)\,=\,\sum_{n\geq 0}a_n^1 Y^nX\cdot X+
\sum_{n\geq 0}b_n^1 Y^n\delta_1\cdot X,$$
and the sums in the previous expression are finite. We already computed this 
kind of sum in the previous Proposition, and in the same way as there, it turns
out that $a_n^1=0$ for all $n\geq 1$ and $b_n^1=0$ for all $n\geq 0$.
With identical reasoning, $d_1^{1,1}c(e_1,e_2,e_5)$ shows that
$a_n^2=0$ for all $n\geq 1$ and $b_n^2=0$ for all $n\geq 0$. Going on like 
this, our cochains looks finally like
$$c\,=\,a_0^1 X\otimes e^2\wedge e^3+a_0^2 X\otimes 
(e^2\wedge e^5 - e^3\wedge e^4)+\ldots.$$ 
Now evaluating $dc$ on $(e_2,e_3,e_4)$ gives
$$dc(e_2,e_3,e_4)\,=\,\pm a_0^1 X\cdot\delta_4 \pm a_0^2 X\cdot\delta_2,$$
and as $c$ should be a cocycle, we conclude $a_0^1=a_0^2=0$, because 
$X\cdot\delta_4=\delta_5$ and $X\cdot\delta_2=\delta_3$. Similarly, we 
can conclude that all $a_0^i$ for $i\geq 1$ are zero by evaluating $dc$ on 
the different $e$-triples. Finally 
$$E_{\infty}^{1,1}\,=\,0.$$

Now let us compute $E_{\infty}^{0,2}$. By Feigin-Fuchs, we have:
$$E_1^{0,2}\,=\,\left(\bigoplus_{n\geq 0}\R Y^nX^2\oplus\bigoplus_{n\geq 0}\R
Y^n\delta_2\oplus\bigoplus_{n\geq 0}\R Y^nX\delta_1\oplus\bigoplus_{n\geq 0}\R 
Y^n\delta_1^2\right)\otimes_{\R}H^2(\fm_0),$$
and a general cochain $c\in E_1^{0,2}$ looks like
\begin{eqnarray*}
c&=&\sum_{n\geq 0}a_n^1 Y^nX^2\otimes e^2\wedge e^3+\sum_{n\geq 0}b_n^1 
Y^n \delta_2\otimes e^2\wedge e^3  \\
&+& \sum_{n\geq 0}c_n^1 Y^nX\delta_1\otimes e^2\wedge e^3+\sum_{n\geq 0}d_n^1 
Y^n \delta_1^2\otimes e^2\wedge e^3  \\
&+&\sum_{n\geq 0}a_n^2 Y^nX^2\otimes (e^2\wedge e^5 - e^3\wedge e^4)+
\sum_{n\geq 0}b_n^1 Y^n \delta_2\otimes (e^2\wedge e^5 - e^3\wedge e^4)  \\
&+& \sum_{n\geq 0}c_n^1 Y^nX\delta_1\otimes (e^2\wedge e^5 - e^3\wedge e^4)+
\sum_{n\geq 0}d_n^1 Y^n \delta_1^2\otimes (e^2\wedge e^5 - e^3\wedge e^4)
+\ldots  
\end{eqnarray*}
We compute once again
$$d_1^{0,2}c(e_1,e_2,e_3)\,=\,\sum_{n\geq 0}a_n^1 Y^nX^2\cdot X+
\sum_{n\geq 0}b_n^1 Y^n \delta_2\cdot X  
+ \sum_{n\geq 0}c_n^1 Y^nX\delta_1\cdot X+\sum_{n\geq 0}d_n^1 
Y^n \delta_1^2\cdot X.$$
It is clear that the first sum cannot mix with the others, as there are no 
$\delta$'s in it. Thus $a_n^1=0$ for all $n\geq 1$. For the other terms, we have
$$Y^n\delta_2\cdot X\,=\,-Y^n\delta_3+Y^{n-1}X\delta_2+\ldots+YXY^{n-2}\delta_2
+XY^{n-1}\delta_2,$$
$$Y^nX\delta_1\cdot X\,=\,-Y^nX\delta_2+Y^{n-1}X^2\delta_1+\ldots+YXY^{n-2}
X\delta_1+XY^{n-1}X\delta_1,$$
and
$$Y^n\delta_1^2\cdot X\,=\,-2Y^n\delta_2\delta_1+Y^{n-1}X\delta_1^2+\ldots+
YXY^{n-2}\delta_1^2+XY^{n-1}\delta_1^2.$$
All sums over $n$ are finite, therefore let us consider only the highest order in 
$Y$. In the fourth sum, there are always two $\delta$s in the highest order
term, they cannot match with the others and therefore $d_n^1=0$ for all 
$n\geq 0$. The terms in the second and third sum can match, but there is then a
term coming from commuting one term into the other, which makes one highest 
coefficient zero. But then the other coefficient must be zero, too. Finally,
$b_n^1=c_n^1=0$ for all $n\geq 0$. The only term which can possibly be non-zero
is thus $a_0^1$. The same reasoning applies to the other $a_n^i$, $b_n^i$,
$c_n^i$ and $d_n^i$ with $i>1$. We are then left with a cochain of the form
$$c\,=\,a_0^1 Y^nX^2\otimes e^2\wedge e^3
+a_0^2 Y^nX^2\otimes (e^2\wedge e^5 - e^3\wedge e^4)+ \ldots .$$
As before, we can then apply $dc$ to other $e_i$-triples in order to show that 
the $a_0^j$ must all be separately zero. Finally
$$E_{\infty}^{0,2}\,=\,0.$$
\fin

\begin{rem}
The result from the previous three Propositions can be interpreted as follows:
The short exact sequence of augmentation
$$0\to ({\mathcal H}^{\rm ad})^+ \to{\mathcal H}^{\rm ad}\to \R\to 0$$
induces a long exact sequence in cohomology, and by the previous results, we 
deduce $H^0(\fm_0;{\mathcal H}^+)=0$, 
$H^1(\fm_0;{\mathcal H}^+)=\R\delta_1\otimes e^1$, and
$H^2(\fm_0;{\mathcal H}^+)=0$. 

As we will see below, it is this term $\delta_1\otimes e^1$ which gives 
rise to the only non-zero term in $HH^2({\mathcal H})$.
\end{rem}

\subsection{Computation of $H^q(\fh;\mathcal{H}^{\rm ad})$}\label{ad}

Here we perform the second step of the outline, i.e. we compute 
$H^q(\fh;\mathcal{H}^{\rm ad})$ for $q=0,1,2$ knowing 
$H^q(\fm_0;\mathcal{H}^{\rm ad})$ for $q=0,1,2$ from the previous subsection,
via the Hochschild-Serre spectral sequence.

Given a Lie algebra $\fg$, a $\fg$-module $L$ and an ideal $\fk$, the spectral 
sequence converges to $H^*(\fg;L)$. As stated before, the $E_2$-term in this 
spectral sequence is
$$E_2^{p,q}\,=\,H^p(\fg\,/\,\fk;H^q(\fk;L)).$$
In our case with $\fk=\fm_0$, $\fg=\fh$ and $L={\mathcal H}^{\rm ad}$,
recall the $2$-dimensional quotient Lie algebra $\fa=\fh\,/\,\fm_0$
generated by $Y$ and $\delta_1$ with the only (non-trivial) relation 
$[Y,\delta_1]=\delta_1$. Then this gives the following spaces:
\begin{eqnarray*}
E_2^{0,0}&=&H^0(\fa;H^0(\fm_0;{\mathcal H}^{\rm ad}))\,=\,\R.  \\
E_2^{1,0}&=&H^1(\fa;H^0(\fm_0;{\mathcal H}^{\rm ad}))\,=\,H^1(\fa;\R). \\
E_2^{0,1}&=&H^0(\fa;H^1(\fm_0;{\mathcal H}^{\rm ad}))\,=\,
(H^1(\fm_0)\oplus \R\delta_1\otimes e^1)^{\fa}.  \\
E_2^{2,0}&=&H^2(\fa;H^0(\fm_0;{\mathcal H}^{\rm ad}))\,=\,H^2(\fa;\R).  \\
E_2^{1,1}&=&H^1(\fa;H^1(\fm_0;{\mathcal H}^{\rm ad}))\,=\,
H^1(\fa;H^1(\fm_0)\oplus \R\delta_1\otimes e^1).   \\
E_2^{0,2}&=&H^0(\fa;H^2(\fm_0;{\mathcal H}^{\rm ad}))\,=\,
H^2(\fm_0)^{\fa}.
\end{eqnarray*}

\begin{rem}
A big difference between the Lie algebras $\fh$ and $\fm_0$ is that $\fh$ has
a so-called {\it Euler element} and $\fm_0$ does not. This means that $\fh$
possesses a basis consisting of eigenvectors with respect to the adjoint action 
of $Y$, while no element acts (in the adjoint action) diagonally on $\fm_0$.
In other words, the grading of $\fh$ is implemented by a grading element $Y$,
while the grading of $\fm_0$ is not given by an inner derivation. 
The existence of an Euler element in $\fh$ implies its existence in $\fa$.

The grading on both Lie algebras induces (second) gradings on cochains spaces. 
We will follow Fuchs' convention \cite{Fuchs} p. 29, and write for a graded Lie 
algebra $\fg$ and a graded $\fg$-module $A$
\begin{equation}      \label{degree}
C^q_{\lambda}(\fg;A)\,=\,\{c\in C^q(\fg;A)\,|\,
\forall\,X_i\in\fg_{\lambda_i}:\,\,\,\,c(X_1,\ldots,X_q)\in 
A_{\lambda_1+\ldots+\lambda_q-\lambda}\,\}.
\end{equation}
By Theorems $1.5.2$ and $1.5.2a$ in \cite{Fuchs}, the subcomplex of 
degree-$0$-cochains for $\fg=\fh$ and a graded module $A$ admitting a basis 
of eigenvectors for the action of $Y$ is homotopy equivalent to the total 
complex. The same is true for $\fa$, but is not for the Lie algebra $\fm_0$.
\end{rem}

\begin{prop}\qquad $\displaystyle 
H^1(\fa;\R)\,=\,\R Y^*\quad\text{and}\quad  H^2(\fa;\R)\,=\,0$.
\end{prop}

\pr The $2$-dimensional Lie algebra $\fa$ admits a grading where $Y$ is the 
grading element of degree $0$ and $\delta_1$ is of degree $1$. This
grading induces a grading on all cochain spaces. By
Theorem $1.5.2$ in \cite{Fuchs}, for cohomology computations one may restrict 
to the degree-$0$-subcomplex. This means the subcomplex given by $\R\subset
C^0(\fa;\R)$, $\R Y^*\subset C^1(\fa;\R)$ and $0\subset C^2(\fa;\R)$.

This implies the proposition.\fin 

\begin{prop}\qquad $\displaystyle
(H^1(\fm_0)\oplus \R\delta_1\otimes e^1)^{\fa}\,=\,\R\delta_1\otimes
e^1$.
\end{prop}

\pr $Y$ and $\delta_1\in\fa$ act trivially on 
$\R\delta_1\otimes e^1$, because the cochain is of degree $0$ with respect to
$Y$ and $e^1\cdot \delta_1\,=\,\pm e^1([\delta_1,-])\,=\,0$. 

On the other hand, $Y$ acts non-trivially on $1\otimes e^1$ and
$1\otimes e^2\in H^1(\fm_0)$, because $[Y,X]=X$ and $X$ corresponds to $e^1$
and $[Y,\delta_2]=2\delta_2$ in $\fh$ (and $\delta_2$ corresponds to $e_2$). 
Thus only $\delta_1\otimes e^1$ is invariant.\fin

\begin{prop}\qquad $\displaystyle H^2(\fm_0)^{\fa}\,=\,0$.
\end{prop}

\pr Observe that $Y$ acts as a grading element on $H^2(\fm_0)$, i.e. every
element of $H^2(\fm_0)$ is $Y$-eigenvector with non-trivial eigenvalue.
It is clear from this and from the explicit description of $H^2(\fm_0)$ in
Theorem \ref{bigraded} that no non-trivial combination of the generators will 
be $Y$-invariant.\fin 

\begin{prop}\qquad $\displaystyle 
H^1(\fa;H^1(\fm_0)\oplus \R\delta_1\otimes e^1)\,=\,\R 
Y^*\otimes(\delta_1\otimes e^1)$.
\end{prop}

\pr Once again, Lie algebra and module admit a basis consisting of eigenvectors
with respect to $Y$. In this situation, as before, one may restrict to 
the subcomplex of cochains of degree $0$, see \cite{Fuchs}, Theorem $1.5.2a$.

Let us compute the degrees of $1\otimes e^1$ and $1\otimes e^2$ using formula 
(\ref{degree}): $1\otimes e^1$ is of degree $1$, because 
$1\cdot e^1(e_1)\in\R\subset({\mathcal H})_0$.
In the same way, $1\otimes e^2$ is of degree $2$ and 
$(1\otimes e^2)\cdot Y=2(1\otimes e^2)$. 

Now it is clear, as $Y^*$ is of degree $0$ and $\delta_1^*$ of degree $1$, 
that there can be built no cochain of degree $0$ from 
${\rm Hom}(\fa;H^1(\fm_0))$. 

On the other hand, $\delta_1\otimes e^1$ is of degree $0$ with respect to $Y$,
and $Y^*\otimes(\delta_1\otimes e^1)$ is obviously a cocycle 
in $C^1(\fa;H^1(\fm_0)\oplus \R\delta_1\otimes e^1)$.\fin

Now the table looks like:

\begin{eqnarray*}
E_2^{0,0}&=&\R  \\
E_2^{1,0}&=&\R Y^* \\
E_2^{0,1}&=&\R\delta_1\otimes e^1  \\
E_2^{2,0}&=&0  \\
E_2^{1,1}&=&\R Y^*\otimes(\delta_1\otimes e^1)   \\
E_2^{0,2}&=&0
\end{eqnarray*}

As a corollary, we have:

\begin{theo}\label{coh1}\qquad $\displaystyle 
H^1(\fh;{\mathcal H}^{\rm ad})\,=\,HH^1({\mathcal H})$\\
is $2$-dimensional, generated in terms of Lie algebra cocycles by $Y^*$
and $\delta_1\otimes e^1$.
\end{theo}

As the differential
$$d^{1,1}_2:E_2^{1,1}\to E_2^{3,0}$$
is zero by dimensional reasons (because $E_2^{3,0}=H^3(\fa;\R)$ is the degree 
$3$ cohomology of a $2$-dimensional Lie algebra), we have as another corollary:

\begin{theo}\label{coh2}\qquad $\displaystyle
H^2(\fh;{\mathcal H}^{\rm ad})\,=\,HH^2({\mathcal H})$\\
is $1$-dimensional, generated in terms of Lie algebra cocycles by 
$X^*\wedge Y^*\otimes \delta_1$.
\end{theo}

\begin{rem}
Observe that all $p$=const columns with $p>3$ on the second page of the 
Hochschild-Serre spectral sequence are zero as $\fa$ is $2$-dimensional.

For example, in order to compute $HH^3({\mathcal H})$, the only space to 
compute would be $H^3(\fm_0;{\mathcal H}^{\rm ad})$ (then one would have to 
take $\fa$-invariants in order to have $E_2^{0,3}$, and afterwards compute 
$d_2^{0,3}$...), because the spaces
$$E_2^{1,2}\,=\,H^1(\fa;H^2(\fm_0;{\mathcal H}^{\rm ad}))\,=\,
H^1(\fa;H^2(\fm_0))$$
and
$$E_2^{2,1}\,=\,H^2(\fa;H^1(\fm_0;{\mathcal H}^{\rm ad}))\,=\,
H^2(\fa;H^1(\fm_0)\oplus\R\delta_1\otimes e^1)$$
are easily seen to be zero by degree arguments as before.
\end{rem}

\begin{rem}
Computing the Hochschild cohomology $HH^*({\mathcal H};\R)$ with trivial 
coefficients is rather easy along the lines of the above computation of 
$HH^2({\mathcal H})$. Indeed, the third step (involving the Feigin-Fuchs 
spectral sequence) is trivial, thus it remains to explore the Hochschild-Serre
spectral sequence. Putting in the previous computations, its second page looks 
in this case (in part) like:
$$
\begin{array}{lll}
\Strut  E_2^{0,3}=H^3(\fm_0)^{\fa} &  
E_2^{1,3}=H^1(\fa;H^3(\fm_0)) &  E_2^{2,3}=H^2(\fa;H^3(\fm_0)) \\
\Strut  E_2^{0,2}=H^2(\fm_0)^{\fa}=0 \quad& 
 E_2^{1,2}=H^1(\fa;H^2(\fm_0)) \quad&  E_2^{2,2}=H^2(\fa;H^2(\fm_0))\\
\Strut  E_2^{0,1}=H^1(\fm_0)^{\fa}=0 &  E_2^{1,1}=H^1(\fa;\R Y^*) & 
 E_2^{2,1}=H^2(\fa;\R Y^*) \\
\Strut  E_2^{0,0}=\R &  E_2^{1,0}=\R Y^* & E_2^{2,0}=H^2(\fa;\R)=0 
\end{array}
$$
By degree arguments, $H^p(\fm_0)^{\fa}=0$ for all $p\geq 1$ and
$H^2(\fa;H^q(\fm_0))=0$ for all $q\geq 2$. Indeed, Theorem
\ref{bigraded} shows that $H^p(\fm_0)$ for $p\geq 1$ is generated by
cocycles which are never  $Y$-invariant. On the other hand, the
cochains computing $H^2(\fa;H^q(\fm_0))$ must be of the form $Y^*\wedge
\delta_1^*\otimes H^q(\fm_0)$, and are therefore of degree $0$ only if
the cochain in $H^q(\fm_0)$ is of degree $1$. This can  onlu happen for
$q=1$ by Theorem \ref{bigraded}. The second page becomes thus:
$$
\begin{array}{lll}
\Strut E_2^{0,3}=0 &  E_2^{1,3}=H^1(\fa;H^3(\fm_0)) &  E_2^{2,3}=0 \\
\Strut E_2^{0,2}=0 \quad &  E_2^{1,2}=H^1(\fa;H^2(\fm_0)) \quad&  E_2^{2,2}=0\\
\Strut  E_2^{0,1}=0 &  E_2^{1,1}=H^1(\fa;\R Y^*) &  E_2^{2,1}=0 \\
\Strut  E_2^{0,0}=\R &  E_2^{1,0}=\R Y^* & E_2^{2,0}=0 
\end{array}
$$
As a consequence, the spectral sequence collapses at the second page, and one 
obtains therefore for all $l\geq 0$
$$HH^l({\mathcal H};\R)\,=\,H^1(\fa;H^{l-1}(\fm_0)).$$
Once again, by degree arguments and the knowledge of $H^*(\fm_0)$ using
Theorem \ref{bigraded}, $H^1(\fa;H^{l-1}(\fm_0))$ is non-zero only for
$l=0,1$, and the result is
$$HH^l({\mathcal H};\R)\,\cong\,\left\{\begin{array}{ccc} \R & {\rm if} & 
l=0,1,2 \\ 0 & {\rm if} & l\geq 3 \end{array}\right. .$$
\end{rem}
\bigbreak

\end{document}